\documentclass[12pt]{article}
\usepackage{amsmath,amsthm,amsfonts}
\theoremstyle{plain}
\newtheorem{theor}{Theorem}
\theoremstyle{plain}
\newtheorem{theor1}{Proposition}
\theoremstyle{remark}
\newtheorem{rem}{Remark}
\theoremstyle{remark}
\newtheorem{rem1}{Example}
\theoremstyle{plain}
\newtheorem{cor}{Corollary}
\normalbaselines
\begin {document}
\vbox {\vspace {6mm}}
\begin {center}
{ \large \bf About interpolation of subspaces\\
of rearrangement invariant spaces \\
generated by Rademacher system\\ [3mm]}
{\bf S.V. Astashkin\\ [5mm]}
\end {center}

\begin {abstract}
The Rademacher series in rearrangement invariant function spaces "closed"
to the space
$ L _ {\infty}$
are considered. In terms of interpolation theory
of operators a correspondence between such spaces and
spaces of coefficients generated by them is stated. It is proved that
this correspondence is one-to-one. Some examples and applications are
presented.
\end {abstract}

\begin {center}
{\large\bf 1. Introduction}
\end {center}

\def\ab {\cal K} \def\bc {{\cal K} (2 ^ k, x; X _ 0, X _ 1)}
\def\cd {\sum _ {k = 1} ^ {\infty} a _ kr _ k}
\def\de {\int _ {0} ^ {u} x ^ * (s) \, ds}
\def\ef {{\cal K} (t, x; L _ {\infty}, G)} \def\fg {{\cal K} (t, a; l _ 1, l _
2)}
\def\gh {(a _ k) _ {k = 1} ^ {\infty}}

\def\ij {(l _ 1, l _ 2) ^ {\cal K}}
\def\jk {(L _ {\infty}, G) ^ {\cal K}}
\def\kl {\varphi _ a}

\def\mn {[t ^ 2]}
\def\no {{\cal K} (t, Ta; L _ {\infty}, G)} \def\op {\cal P}
\def\pq {\gamma _ f}
\def\qr {\delta _ f}
\def\rv {L _ {\infty}}
\def\vw {{\cal K} (t, x; L _ 1, L _ 2)} \def\wx {{\cal K} (2 ^ k, x; L _
{\infty}, G)}
\def\xy {M (\varphi)}
\def\yz {l _ 1 (\log)}

\def\us {\varphi _ 0}
\def\st {\varphi _ 1}
\def\td {\Lambda _ p (\varphi)}
\def\dn {(X _ 0, X _ 1) ^ {\cal K}}  \def\an {{\cal K} (2 ^ k,
a; l _ 1, l _ 2)}
\def\ao {\vec {X}}
\def\eu {\vec {Y}}
\def\ej {\tilde {\varphi}}

Let
$$
r_ k (t) = \, \rm sign\sin {2 ^ {k-1} {\pi} t} \; \;\; (k = 1,2..) $$
be the Rademacher functions on the segment
$ [0,1]. $
Define the linear operator
$$
Ta (t) = \, \sum _ {k = 1} ^ {\infty} \, a _ kr _ k (t) \; \;\mbox{for}\;\;a = (a _ k) _ {k = 1} ^ {\infty} \in {l _ 2}.\eqno {(1)} $$
It is well-known (see, for example, [1, p.340-342]) that by fixed
$ a $
$Ta(t)$
is an almost everywhere finite function on
$[0,1].$
Moreover,
from Khintchine's inequality it follows that
$$
|| Ta ||_{L_p} \,\asymp {\, || a || _ 2}\;\;\;\mbox{for}\;\;\;1\le p<\infty,\eqno {(2)} $$
where
$\|a\|_p=\,(\sum_{k=1}^\infty|a_k|^p)^{1/p},$
as usual.
The last means the existence of two-sided estimates with constants depending
only on
$p.$
Also, it can easily be checked that
$$
|| Ta || _ {L_\infty} \, = \, || a || _ 1.\eqno {(3)} $$

A more detailed information on the behaviour of Rademacher series
can be obtained by treating them in the framework
of general rearrangement invariant spaces.

Recall that a Banach space
$ X $
of measurable functions
$x=x(t)$
on
$ [0,1] $
is said to be
an rearrangement invariant space (r.i.s.) if the inequality
$ x ^ * (t) \le {y ^ * (t)}$
for
$t\in [0,1]$
and
$ y\in X $
imply
$ x\in X $
and
$ || x ||\le {|| y ||}.$
Here and next
$ z ^ * (t) $
is the nonincreasing rearrangement of a function
$ | z (t) | $
concerning to Lebesgue measure denoted by meas
[2, p.83].

Important examples of r.i.s.'s are Marcinkiewicz and Orlicz spaces.

By
$ \op $
will be denoted the cone of nonnegative increasing concave
functions on the semiaxis
$(0, \infty).$

If
$\varphi\in\op,$
then the Marcinkiewicz space
$M(\varphi)$
consists of all measurable functions
$x=x(t)$
such that
$$
||x||_{M(\varphi)}\,=\,\sup\left\{\frac{1}{\varphi(t)}\,\int_o^t\,x^*(s)\,ds:\,0<t\le 1\right\}\,<\,\infty.$$

If
$S(t)$
is a nonnegative convex continuous function on
$[0,\infty),$
$S(0)=0,$
then the Orlicz space
$L_S$
consists of all measurable functions
$x=x(t)$
such that
$$
||x||_S\,=\,\inf\left\{u>0:\,\int_{0}^{1}\,S(\frac{|x(t)|}{u})\,dt\,\le{\,1}\right\}\,<\,\infty.$$
In particular, if
$ S (t) = t ^ p $
$(1\le p < \infty), $
then
$ L _ S = L _ p $.

For any r.i.s.
$ X $
on
$ [0,1] $
we have
$ L _ {\infty} \subset {X} \subset {L _ 1} $
[2, p.124].
By
$X^0$
will be denoted the closure of
$L_\infty$
in an r.i.s.
$X.$

In problems discussed below, a special role is played by the Orlicz space
$ L _ N, $
where
$ N (t) = \exp (t ^ 2) -1 $
or, more precisely, by the space
$ G=\,L _ N^0.$
In the paper [3], V.A.Rodin and E.M.Semenov proved a theorem about
the equivalence of Rademacher system to the standard basis in the space
$l_2.$

\noindent {\bf Theorem A.} {\it Suppose that
$ X $
is an r.i.s. Then
$$
|| Ta || _ X \; = \; \Bigl\|\cd \Bigr\| _ X \,\asymp {\, || a || _ 2} $$
if and only if
$ X\supset G. $ }

By Theorem A, the space
$ G $
is the minimal space among r.i.s.'s
$X$
such that the Rademacher system is equivalent in
$X$
to the standard basis of
$ l _ 2. $

In this paper, we consider problems related to the behaviour of
Rademacher series in r.i.s.'s intermediate between
$ L _ {\infty} $
and
$ G. $
The main role is played here by concepts and methods of interpolation theory
of operators.

For a Banach couple
$ (X _ 0, X _ 1), x\in {X _ 0 + X _ 1}, $
and
$ t> 0 $
we introduce the Peetre
$ {\cal K} $
-functional
$$
{ \cal K} (t, x; X _ 0, X _ 1) \, = \, \inf \{|| x _ 0 || _ {X _ 0} + t || x _ 1 || _ {X _ 1}: \, x = x _ 0 + x _ 1, x _ 0\in {X _ 0},x_1\in X_1 \}. $$

Let
$ Y _ 0 $
be a subspace of
$ X _ 0 $
and
$ Y _ 1 $
a subspace of
$ X _ 1 $.
A couple
$ (Y _ 0, Y _ 1) $
is called a
$ {\cal K} $
-subcouple of a couple
$ (X _ 0, X _ 1) $
if
$$
{ \cal K} (t, y; Y _ o, Y _ 1) \, \asymp{{\,\cal K} (t, y; X _ 0, X _ 1)}, $$
with constants independent of
$ y\in {Y _ 0 + Y _ 1} $
and
$ t> 0. $

In particular, if
$ Y _ i \, = \, P (X _ i), $
where
$ P $
is a linear projector bounded from
$ X _ i $
into itself for
$i=0,1,$
then
$ (Y _ 0, Y _ 1) $
is a $ {\cal K} $
-subcouple of
$ (X _ 0, X _ 1) $
(see [5] or [6, p.136]). At the same time,
there are many examples of subcouples that are not
$ {\cal K} $
-subcouples (~see [6, p.589], [7],
and Remark 2 of this paper).

We shall consider the case:
$ X _ 0 \, = L _ {\infty}, \, X _ 1 = \, G, \, Y _ 0 = \, T (l _ 1),$
and
$ Y _ 1 = \, T (l _ 2), $
where
$ T $
is given by (1). From (3) and Theorem A it follows that
$$
{ \cal K} (t, Ta; T (l _ 1), T (l _ 2)) \, \asymp {\, {\cal K} (t, a; l _ 1, l _ 2)}. \eqno {(4)} $$
In spite of the fact that
$ T (l _ 1) $
is uncomplemented in
$ L _ {\infty} $
(see [8] or [9, p.134]) the following assertion holds.

\begin{theor}
The couple
$ (T (l _ 1), T (l _ 2)) $
is a
$ {\cal K} $
-subcouple of the couple
$ (L _ {\infty}, G). $
In other words (see (4)),
$$
{ \cal K} (t, Ta; L _ {\infty}, G) \, \asymp {{\cal K} (t, a; l _ 1, l _ 2)},\eqno {(5)} $$
with constants independent of
$ a = (a _ k) _ {k = 1} ^{\infty} \in {l _ 2} $
and
$ t> 0. $
\end{theor}

We shall use in the proof of Theorem 1 an assertion about the
distribution of Rademacher sums. It was proved by S.Montgomery-Smith [10].

\noindent {\bf Theorem B.} {\it There exists a constant
$ A\ge 1 $
such that for all
$ a = \gh\in {l _ 2} $
and
$ t> 0 $
$$
{\rm meas}\Bigl \{s\in {[0,1]}: \,\cd (s) > A ^ {-1} \kl (t) \Bigr\} \,\ge {A ^{-1} \exp (-At ^ 2)}, \eqno {(6)} $$
where
$ \kl (t) \, = \, \fg. $}
\vskip 0.2cm

We need now some definitions from interpolation theory of operators.

We say that a linear operator
$ U $
is bounded from a Banach couple
$ \ao = (X _0, X _ 1) $
into a Banach couple
$ \eu = (Y _ 0, Y _ 1) $
(in short,
$U: \,\ao\to {\eu})$
if
$U$
is defined on
$X_0+X_1$
and acts as bounded operator from
$ X _ i $
into
$ Y _ i $
for
$ i=0,1.$

Let
$ \ao = (X _ 0, X _ 1) $
be a Banach couple. A space
$ X $
such that
$ X _ 0\cap{X _ 1} \subset X\subset {X _ 0 + X _ 1} $
is called an interpolation space between
$ X _ 0 $
and
$ X _ 1 $
if each linear operator
$ U: \,\ao\to {\ao} $
is bounded from
$ X$
into itself.

For every
$p\in [1,\infty],$
we shall denote by
$ l _ p (u _ k), \, u _ k\ge 0 \, (k = 0,1..) $
the space of all two-sided sequences of real numbers
$ a = (a _ k) _ {k = -\infty} ^ {\infty} $
such that the norm
$|| a || _ {l _ p (u _ k)} \, = \, || (a _ ku _ k) || _ p $
is finite.
Let
$ E $
be a Banach lattice of two-sided sequences,
$ (\min(1,2^k))_{k=-\infty}^\infty\in E. $
If
$ (X _ 0, X _ 1)$
is a Banach couple, then the space of the real
$ \ab $
-method of interpolation
$\dn_E$
consists of all
$ x\in {X _ 0 + X _ 1} $
such that
$$
|| x || \, = \, || (\bc) _ k || _ E\,<\,\infty.$$
It is readily checked that the space
$ \dn _ E $
is an interpolation space between
$ X_ 0 $
and
$ X _ 1 $
(see, for example, [11, p.422]). In the special case
$ E= \, l _ p (2 ^ {-k\theta}) $
$ (0 < \theta < 1,$
$1\le p\le {\infty}) $
we obtain the spaces
$ (X _ 0, X _ 1) _ {\theta, p} $
(for the detailed exposition of their properties see [4]).

A couple
$ \ao = (X _ 0, X _ 1) $
is said to be a
$ \ab $
-monotone couple if
for every
$ x\in {X _ 0 + X _ 1} $
and
$ y\in {X _ 0 + X _ 1} $
there exists a linear operator
$ U: \,\ao\to {\ao} $
such that
$ y = Ux$
whenever
$$
{\ab} (t, y; X _ 0, X _ 1) \, \le {\, {\ab} (t, x; X _ 0, X _ 1)}\;\;\mbox{for all}\;\;t> 0.$$

As it is well-known (see, for example, [11, p.482]), any interpolation space
$ X$
with respect to a
$ \ab $
-monotone couple
$ (X _ 0, X _ 1)$
is described by the real
$\ab$
-method. It means that for some
$E$
$$
X \, = \,\dn _ E. $$

In particular, by the Sparr theorem [12] the couple
$ (l _ 1, l _ 2) $
is a
$ \ab $
-monotone couple. Therefore, if
$ F $
is an interpolation space between
$l_1$
and
$l_2,$
then there exists
$E$
such that
$$
F \, = \,\ij _ E. $$
Hence Theorem 1 allows to find an r.i.s. that
contains Rademacher series with coefficients belonging
to an arbitrary interpolation space between
$ l _ 1 $
and
$ l _ 2.$
In the paper [3], the similar result was obtained for sequence spaces
satiafying more restrictive conditions (see Remark 3).

\begin{theor}
Let
$ F $
be an interpolation sequence space between
$ l _ 1 $
and
$ l _ 2$
and
$ F = \, \ij _ E. $
Then for the r.i.s.
$ X = \, \jk _ E $
we have
$$
\Bigl\|\cd\Bigr\| _ X \,\asymp {|| a || _ F}, $$
with constants independent of
$ a = \gh. $
\end{theor}

The last result shows that a correspondence given by Theorem 2 is one-to-one.

\begin{theor}
Let r.i.s.'s
$ X _ 0 $
and
$ X _ 1 $
be two interpolation spaces between
$ \rv $
and
$ G $.
If
$$
\biggl\|\cd\biggr\| _ {X _ 0} \, \asymp {\,\biggl\|\cd\biggr\| _ {X _ 1}}, $$
then
$ X _ 0 \, = \, X _ 1$
and the norms of
$X_0$
and
$X_1$
are equivalent.
\end{theor}

In the papers [13] and [13], the similar results were obtained by
additional conditions with respect to spaces
$X_0$
and
$X_1.$
\vskip 0.6cm

\begin {center}
{\large \bf 2. Proofs}
\end {center}
\begin{proof}[Proof of Theorem 1.]
It is known [2, p.164] that the
$ \ab $
-functional of a couple of Marcinkiewicz spaces is given by the formula
$$
{\ab} (t, x; M (\us), M (\st)) \, = \, \sup _ {0 < u\le 1} {\frac {\de} {\max(\us (u), \st (u) /t)}}. $$
If
$ N (t) = \exp (t ^ 2) -1, $
then the Orlicz space
$ L _ N $
coincides with the Marcinkiewicz space
$ M (\varphi _ 1),$
where
$\varphi_ 1 (u) = u {\log} _ 2 ^ {1/2} (2/u) $
[3]. In addition,
$ \rv = M (\varphi _ 0)$
where
$\varphi _ 0 (u) = u.$
Therefore,
$$ \ef \, = \,\sup _ {0 < u\le 1} \left \{\frac{1}{u}\,\de \,\min (1, t {\log} _ 2 ^{-1/2} (2/u)) \right\}\;\;\mbox{for}\;\;x\in G.\eqno {(7)} $$
Since
$ x ^ * (u) \le {1/u\de}, $
then from (7) it follows
$$ \ef \,\ge {\, \sup _ {k = 0,1..} \, \left \{x ^ * (2 ^ {-k}) \min (1, t (k+ 1) ^ {-1/2}) \right \}}. $$
Hence,
$$
\ef \,\ge {\, x ^ * (2 ^ {-k _ t})}\;\;\mbox{for}\;\;t\ge 1, \eqno {(8)} $$
where
$ k _ t = \mn-1 $
$([z] $
is the integral part of a number
$ z $).

Let now
$a = \gh\in\,l _ 2$
and
$ x (t) = Ta (t) = \cd (t).$
By the Holmstedt formula [14],
$$ \kl (t) \, \le {\sum _ {k = 1} ^ {[t ^ 2]} a _ k ^ * \, + \, t\left \{\sum_ {k = \mn + 1} ^ {\infty} (a _ k ^ *) ^ 2\right \} ^ {1/2}} \, \le {\, B\kl(t)}, \eqno {(9)} $$
where
$ \kl (t) = \, \fg,$
$ (a _ k ^ *)_{k=1}^\infty$
is nonincreasing rearrangement of the sequence
$ (| a _ k |) _ {k = 1} ^{\infty}, $
and
$B>0$
is a constant independent of
$a = \gh$
and
$t>0.$

Assume, at first, that
$ a\not\in {l _ 1}. $
Then inequality (9) shows that
$$
\lim _ {t\to {0 +}} \kl (t) \, = \, 0 \mbox{\;\;\;and\;\;\;}\lim _ {t\to {\infty}} \kl (t) \,= \, \infty. $$
The function
$ \kl$
belongs to the class
$\op$
[4, p.55]. Therefore it
maps the semiaxis
$ (0, \infty) $
onto
$ (0, \infty) $
one-to-one, and there exists the inverse function
$ \kl ^ {-1} $.
By Theorem B, we have
$$ n _ {| x |} (\tau) \, = \, \rm meas \{s\in {[0,1]}: \, | x (s) | > \tau \} \,\ge{\, \psi (\tau)}\;\;\mbox{for}\;\;\tau > 0, $$
where
$ \psi (\tau) = A ^ {-1} \, \exp \{-A[\kl ^ {-1} (\tau {A})] ^ 2 \}.$
Passing to rearrangements we obtain
$$
x ^ * (s) \, \ge {\, \psi ^ {-1} (s)}\;\;\;\;\;\mbox{for}\;\;0 < s < A ^ {-1}. \eqno {(10)} $$
Obviously, by condition
$ t\ge {C _ 1} = C _ 1 (A) =\sqrt {2\log _ 2 (2A)}, $
it holds
$$
2 ^ {-k _ t/2} \, < \, A ^ {-1}, \eqno {(11)} $$
for
$ k _ t = \mn-1. $

Hence (8) and (10) imply
$$
\ef \,\ge {\, \psi ^ {-1} (2 ^ {-k _ t})}. \eqno {(12)} $$
Combining the definition of the function
$ \psi$
with (11), we obtain
\begin{multline*}
\psi ^ {-1} (2 ^ {-k _ t}) = A ^ {-1} \kl\left (A ^ {-1/2} \ln ^ {1/2} (A ^{-1} 2 ^ {k _ t}) \right) \ge {A ^ {-1} \kl\left (\sqrt {k _ t\ln 2 / (2A)}\right)} \ge\\
\ge {A ^ {-3/2} \sqrt {\ln 2/2} \kl (\sqrt {k _ t})} \ge {A ^ {-3/2} \sqrt {\ln2/2} t ^ {-1} \sqrt {k _ t} \kl (t)}.
\end{multline*}
From the inequality
$ t\ge {C _ 1} \ge {\sqrt {2}} $
it follows
$$ \frac {\sqrt {k _ t}} {t} \, \ge {\, \frac {\sqrt {\mn-1}} {\sqrt {\mn +1}}} \, \ge {\, 3 ^ {-1/2}}. $$
Therefore, by (12), we have
$$
\ef \,\ge {\, C _ 2\kl (t)}\;\;\;\mbox{for}\;\;t\ge {C _ 1},$$
where
$ C _ 2 = C _ 2 (A) = \sqrt {\ln 2/6} \, A ^ {-3/2}. $

If now
$ t\ge 1, $
then the concavity of the
$ \ab $
-functional
and the previous inequality yield
$$
\ef \,\ge {\, C _ 1 ^ {-1} {\ab} (tC _ 1, x; \rv, G)} \, \ge {\,\frac{C _ 2}{C _ 1}{\kl} (C _ 1t)} \, \ge {\,\frac{C _ 2}{C _ 1} {\kl} (t)}. $$
Using the inequalities
$ || a || _ 2\le {|| a || _ 1} \, (a\in {l _ 1}) $
and
$ || x || _ G\le {|| x|| _ {\infty}}\, (x\in {\rv}), $
the definition of the
$ \ab $
-functional, and Theorem A, we obtain
$$
\ef \, = \, t || x || _ G \,\ge {\, C _ 3t || a || _ 2} \, = \, C _ 3\kl (t)\;\;\;\mbox{for}\;\;0<t\le 1.$$
Thus,
$$
\fg \,\le {\, C\no}, \eqno {(13)} $$
if
$ C = \max (C _ 3 ^ {-1}, C _ 1/C _ 2). $

Suppose now
$ a\in {l _ 1}. $
By (9), without loss of generality, we can assume that the function
$ \kl $
maps the semiaxis $ (0, \infty) $
onto the interval
$ (0, || a || _1) $
one-to-one. Hence we can define the mappings
$ \kl ^ {-1}: \, (0, || a || _ 1) \to {(0,\infty)},$
$\psi: \, (0, A ^ {-1} || a || _ 1) \to {(0, A ^ {-1})}, $
and
$\psi ^ {-1}: \, (0, A ^ {-1}) \to {(0, A ^ {-1} || a || _ 1)}. $
Arguing as above, we get inequality (13).

The opposite inequality follows from Theorem A and relation (3).
Indeed,
$$
\no \,\le {\, \inf \{|| Ta ^ 0 || _ {\infty} + t || Ta ^ 1 || _ G: \, a = a^ 0 + a ^ 1, a ^ 0\in {l _ 1}, a ^ 1\in {l _ 2} \}} \le $$
$$
\le {D\fg}. $$
\end{proof}
\begin{proof}[Proof of Theorem 2.]
 It is sufficient to use Theorem 1 and the definition of the real
 $ \ab $
-method of interpolation.
 \end{proof}
\vskip 0.3cm

For the proof of Theorem 3 we need some definitions and auxiliary
assertions. These results also are of some independent interest.

Let
$f(t)$
be a fuction defined on the interval
$(0,l),$
where
$l=1$
or
$l=\infty.$
Then the dilation function of
$f$
is defined as follows:
$$
{ \cal M} _ f (t) \, = \, \sup\left \{\frac {f (st)} {f (s)}: \, s, st\in {(0,l)} \right \}, \; \;\;\mbox{if}\;\; t\in {(0, l)}. $$
Since this function is semimultiplicative, then
there exist numbers
$$
\pq \, = \,\lim _ {t\to {0 +}} \frac {\ln {{\cal M} _ f (t)}} {\ln t} \mbox{      and     }\qr\, = \,\lim _ {t\to {\infty}} \frac {\ln {{\cal M} _ f (t)}} {\ln t}. $$

A Banach couple
$ \ao = (X _ 0, X _ 1) $
is called a partial retract of a couple
$ \eu = (Y _ 0, Y _ 1) $
if each element
$ x\in {X _ 0 + X _ 1} $
is orbitally equivalent to some element
$ y\in {Y _ 0 + Y _ 1}. $
The last means that there exist linear operators
$ U: \,\ao\to {\eu} $
and
$ V: \,\eu\to {\ao}$
such that
$ Ux = y $
and
$ Vy = x. $

\begin{theor1}
Suppose
$ \xy $
is an Marcinkiewicz space on
$ [0,1]. $
If
$ \gamma_\varphi > 0,$
then
$ \ao = (\rv, \xy) $
is a
$ \ab $
-monotone couple.
\end{theor1}
\begin{proof}
It is sufficient to show that the couple
$ \ao $
is a partial retract of the couple
$ \eu = (\rv, \rv (\ej)),$
where
$$
|| x || _ {\rv (\ej)} \, = \, \sup _ {0 < t\le 1} \ej (t) | x (t) |\;,\;\ej (t)= t/\varphi (t). $$
Indeed,
a partial retract of a
$ \ab $
-monotone couple is a
$ \ab $
-monotone couple [11, p.420], and by the Sparr theorem [12],
$ \eu $
is a
$ \ab $
-monotone couple.

First note that
the inclusion
$ \rv\subset \xy $
implies
$ \rv + \xy = \xy $.
So, let
$ x\in {\xy}. $
Without loss of generality [2, p.87], assume that
$x (t) = x ^ * (t). $
Define the operator
$$
U _ 1y (t) = \, \sum _ {k = 1} ^ {\infty} 2 ^ k\int _ 0 ^ {2 ^ {-k}} y (s)ds \,\chi _ {(2 ^ {-k}, 2 ^ {-k + 1}]} (t)\;\;\;\mbox{for}\;\;y\in\xy. $$
The concavity of the function
$ \varphi$
and properties of the nonincreasing rearrangement imply
$$
|| U _ 1y || _ {\rv (\ej)} \, \le {\, 2\sup _ {k = 1,2..} (\varphi (2 ^ {-k +1})) ^ {-1} \int _ 0 ^ {2 ^ {-k}} y ^ * (s) ds} \, \le {\, 2 || y || _ {\xy}},$$
Hence
$ U _ 1: \,\ao\to {\eu}. $
Since
$ x (t) $
nonincreases, then
$ U _ 1x (t) \ge {x (t)}. $
Therefore the linear operator
$$
Uy (t) = \, \frac {x (t)} {U _ 1x (t)} U _ 1y (t) $$
is bounded from the couple
$ \ao $
into the couple
$ \eu.$
In addition,
$ Ux (t) = x (t). $

Take for
$ V $
the identity mapping, i.e.,
$ Vy (t) = y (t). $
Since
$ \pq > 0, $
then, by [2, p.156], we have
$$
|| Vy || _ {\xy} \, \le {\, C\sup _ {0 < t\le 1} \ej (t) y ^ * (t)} \, \le {\,C\sup _ {0 < t\le 1} \ej (t) | y (t) |} \, = \, C || y || _ {\rv (\ej)}. $$
Therefore
$ V: \,\eu\to {\ao} $
and
$ Vx = x. $

Thus an arbitrary element
$ x\in {\xy} $
is orbitally equivalent to itself as to element of the space
$ \rv + \rv (\ej). $
This completes the proof.
\end{proof}

\begin{cor}
If
$\gamma_\varphi>0,$
then
$(L_\infty,M(\varphi)^0)$
is a
$ \ab $
-monotone couple.
\end{cor}
\begin{proof}
Assume that
$x$
and
$y$
belong to the space
$\xy ^ 0$
and
$$
{\ab} (t, y; \rv, \xy ^ 0) \, \le {\, {\ab} (t, x; \rv, \xy ^ 0)} \;\;\;\mbox{for}\;\;t> 0. $$
If
$ z\in {\xy ^ 0},$
then
$$
{\ab} (t, z; \rv, \xy ^ 0) \, = \, {\ab} (t, z; \rv, \xy).$$
Therefore,
$$
{\ab} (t, y; \rv, \xy) \, \le {{\ab} (t, x; \rv, \xy)}\;\;\;\mbox{for}\;\;t>0.$$
Hence, by Proposition 1, there exists an operator
$ T: \, (\rv, \xy) \to{(\rv, \xy)} $
such that
$ y = Tx. $
It is readily seen that
$ \xy ^ 0 $
is an interpolation space concerning to the couple
$ (\rv, \xy). $
Therefore
$ T: \, (\rv, \xy ^ 0) \to {(\rv, \xy ^ 0)}.$
\end{proof}

We define now two subcones of the cone
$ \op. $
Let us denote by
${\op} _ 0 $
the set of all functions
$ f\in \op$
such that
$ \lim _ {t\to {0 +}} f (t) \, = \, \lim _ {t\to {\infty}} f (t) /t\, = \, 0. $
If
$ f\in {\op}, $
then
$ 0\le {\pq} \le {\qr} \le 1 $
[2, p.76]. Let
$ {\op} ^ {+ -} $
be the set of all
$f\in {\op} $
such that
$ 0 < \pq\le {\qr} < 1. $
It is obvious that
$ {\op} ^ {+-} \subset {{\op} _ 0}. $

A couple
$ (X _ 0, X _ 1) $
is called a
$ {\ab} _ 0 $
-complete couple if for any function
$ f\in {{\op} _ 0} $
there exists an element
$ x\in {X _ 0 + X _ 1} $
such that
$$
{\ab} (t, x; X _ 0, X _ 1) \, \asymp {\, f (t)}. $$
In other words, the set
$ {\ab} (X _ 0 + X _ 1) $
of all
$ \ab $
-functionals of a
$ {\ab} _0 $
-complete couple
$ (X _ 0, X _ 1) $
contains,
up to equivalence,
the whole of the subcone
$ {\op} _ 0.$

\begin{theor1}
The Banach couple
$ (L _ 1 (0, \infty), L _ 2 (0, \infty)) $
is a
$ {\ab} _ 0 $
-complete couple.
\end{theor1}
\begin{proof}
By the Holmstedt formula for functional spaces [14],
$$ \vw \,\asymp {\, \max\left \{\int _ 0 ^ {t ^ 2} x ^ * (s) \, ds, t\left[\int _ {t ^ 2} ^ {\infty} (x ^ * (s)) ^ 2 \, ds\right] ^ {1/2} \right \}}\eqno{(14)} $$
If
$ f\in {{\op} _ 0}, $
then
$ g (t) = \, f (t ^ {1/2}) $
belongs to
$ {\op} _ 0,$
also. Let us denote
$ x (t) = \, g ' (t).$
Then
$ x (t) = x ^ * (t)$
and
$$
\int _ 0 ^ tx (s) \, ds \, = \, g (t). \eqno {(15)} $$
Assume that
$ f\in {\op ^ {+ -}}.$
If
$\qr < 1,$
then there exists
$ \varepsilon > 0 $
such that for some
$ C> 0$
$$
G (s) \, = \, f (s ^ {1/2}) \, \le {\, C (\sqrt {s/t}) ^ {1-\varepsilon} \, f (t ^{1/2})}, \; \;\;\mbox{if}\;\;\;s\ge t. $$
Since
$ g\in {{\op} _ 0}, $
then
$ g ' (t) \le {g (t)/t}. $
Therefore for
$t>0$
$$ \int _ t ^ {\infty} (x (s)) ^ 2 \, ds \,\le {\, \int _ t ^ {\infty} \frac {g^ 2 (s)} {s ^ 2} \, ds} \, \le {\, C ^ 2t ^ {\varepsilon-1} (f (t ^ {1/2})) ^ 2\,\int _ t ^ {\infty} s ^ {-1-\varepsilon} \, ds} \, = \, C ^ 2 {\varepsilon \, t} ^{-1} (g (t)) ^ 2.$$
Combining this with (14) and (15) we obtain
$$
\vw \,\asymp {\, g (t ^ 2)} \, = \, f (t).$$
Thus
$ {\ab} (L _ 1 + L _ 2) \supset {{\op} ^ {+ -}}. $
Hence, in particular, the intersection
$ {\ab} (X _ 0 + X _ 1) \cap {\op ^ {+ -}} $
is not empty. Therefore, by [15, 4.5.7],
$(L_1,L_2)$
is a
${\ab}_0$
-complete Banach couple. This completes the proof.
\end{proof}

Let
$ {\ab} (l _ 1 + l _ 2) $
be the set of all
$ \ab $
-functionals corresponding to the couple
$ (l _ 1, l _ 2).$
By
$ {\cal F} $
we denote the set of all functions
$ f\in {\op} $
such that
$$
f (t) = f (1) t \;\;\;\mbox{for}\;\;0 < t\le 1\;\;\;\;\mbox{and}\;\;\;\; \lim _ {t\to {\infty}} f (t) /t = 0.$$

\begin{cor}
Up to equivalence,
$$
{\ab} (l _ 1 + l _ 2) \supset {{\cal F}}. $$
\end{cor}
\begin{proof}
It is well-known (see, for example, [4, p.142])
that for
$x\in {L _ 1 (0, \infty) + L _ {\infty} (0, \infty)} $
and
$u>0$
$$
{\ab} (u, x; L _ 1, L _ {\infty}) \, = \, \de.\eqno {(16)} $$
In addition,
$$ L _ 1 \, = \, (L _ 1, L _ {\infty}) _ {l _ {\infty}} ^ {\ab}\;\;\;\;\mbox{and}\;\;\;\;L _ 2 \, =\,(L _ 1, L _ {\infty}) _ {l _ 2 (2 ^{-k/2})} ^ {\ab}. $$

The spaces
$ l _ {\infty} $
and
$ l _ 2 (2 ^ {-k/2}) $
are interpolation spaces concerning to the couple
$ \vec {l _ {\infty}} = (l _ {\infty}, l _ {\infty} (2 ^{-k})) $
[4]. Therefore, by the reiteration theorem (see [16] or [17]),
$$
\vw \,\asymp {{\ab} (t, {\ab} (\cdot, x; L _ 1, L _ {\infty}); l _{\infty}, l _ 2 (2 ^ {-k/2}))}\;\;\;\mbox{for}\;\;x\in {L _ 1 + L _ 2}\eqno {(17)} $$

Introduce the average operator:
$$
Qx (t) = \, \sum _ {k = 1} ^ {\infty} \, \int _ {k-1} ^ kx (s) ds \,\chi _{(k-1, k]} (t), \;\;\;\mbox{if}\;\;\;t> 0. $$
From (16) it follows that
$$
{\ab} (t, Qx ^ *; L _ 1, L _ {\infty}) \, = \, {\ab} (t, x; L _ 1, L _{\infty}), $$
for all positive integer
$t.$
Both functions in the last equation are concave. Therefore,
$$
{\ab} (t, Qx ^ *; L _ 1, L _ {\infty}) \, \asymp {\, {\ab} (t, x; L _ 1. L _{\infty})} \; \;\;\;\mbox{for all}\;\;\;t\ge 1. $$
Hence (17) yields
$$
{\ab} (t, Qx ^ *; L _ 1, L _ 2) \, \asymp {\, {\ab} (t, x; L _ 1, L _ 2)}, \; \;\;\;{if}\;\;\;t\ge 1. \eqno {(18)} $$

Let now
$ f\in {\cal F}. $
Since
$ {\cal F} \subset {{\op} _ 0}, $
then, by Proposition 2, there exists a function
$ x\in {L _ 1 (0, \infty) + L _2 (0, \infty)} $
such that
$$
{\ab} (t, x; L _ 1, L _ 2) \, \asymp {\, f (t)}. \eqno {(19)}$$

Clearly, the operator
$Q$
is a projector in the spaces
$ L _ 1 $
and
$ L _ 2 $
with norm 1. Moreover,
$ Q (L _ 1) = l _1$
and
$Q (L _ 2) = l _ 2.$
Hence, by the theorem about complemented subcouples mentioned in
Introduction (see [5] or [6, p.136]),
$$
{\ab} (t, Qx ^ *; L _ 1, L _ 2) \, \asymp {\, \fg}\;\;\;\;\mbox{for}\;\;\;t> 0,$$
where
$ a = \, (\int _ {k-1} ^ kx ^ * (s) ds) _ {k = 1} ^ {\infty}. $

Thus (18) and (19) imply
$$
{\ab} (t, a; l _ 1, l _ 2) \, \asymp {\, f (t)}\;\;\;\;\mbox{for}\;\;\;t\ge 1. $$
The last relation also holds if
$ 0 < t\le 1.$
Indeed, in this case
$$
\fg = t || a || _ 2 = t {\ab} (1, a; l _ 1, l _ 2) \asymp {tf (1)} = f (t). $$
This completes the proof.
\end{proof}

\begin{proof}[Proof of Theorem 3.]
As it was already mentioned in the proof of Theorem 1, the Orlicz space
$ L_ N,$
$N (t) = \exp (t ^ 2) -1,$
coincides with the Marcinkiewicz space
$ M(\varphi _ 1),$
for
$ \varphi _ 1 (u) = \, u\log _ 2 ^ {1/2} (2/u). $
Since
$ \gamma_ {\varphi _ 1} = 1,$
then Corollary 1 implies that the couple
$(\rv, G) $
is a
$ \ab $
-monotone couple. Hence,
$$
X _ 0 = \, \jk _ {E _ 0}\;\;\;\;\mbox{and}\;\;\;\; X _ 1 = \, \jk _ {E _ 1}, \eqno {(20)} $$
for some parameters of the real
$ \ab $
-method of interpolation
$ E _ 0 $
and
$ E _ 1. $
By Theorem 2,
$$
\biggl\|\cd \biggr\| _ {X _ i} \, \asymp {\, || (a _ k) || _ {F _ i}}, $$
where
$ F _ i = \ij _ {E _ i} \; (i = 0,1).$
So, by condition,
$$
\ij _ {E _ 0} \, = \, \ij _ {E _ 1}. \eqno {(21)} $$
The last means that the norms of spaces
$ E _ 0 $
and
$ E _ 1 $
are equivalent on the set
$ {\ab} (l _ 1 + l _2).$
It is readily to check that
this set coincides, up to the equivalence, with the set
$ {\ab} (\rv + G)$
of all
$ \ab $
-functionals corresponding to the couple
$ (\rv, G). $
More precisely,
$$
{\ab} (l _ 1 + l _ 2) \, = \, {\ab} (\rv + G) \, = \, {\cal F}. \eqno {(22)} $$
In fact, by Theorem 1 and Corollary 1,
$ {\cal F} \subset {{\ab} (l _ 1 + l _ 2)} \subset {{\ab} (\rv + G)}. $
On the other hand, since
$ \rv\subset G $
with the constant 1 and
$ \rv $
is dense in
$ G, $
then
$ {\ab} (\rv + G) \subset {\cal F} $
[11, p.386].

Let now
$ x\in {X _ 0}. $
By (20), we have
$ (\wx) _ k\in {X _ 0}. $
Using (22), we can find
$ a\in {l _2} $
such that
$$
\an \,\asymp {\, \wx},$$
for all positive integer
$k.$
Since a parameter of
$ \ab $
-method is a Banach lattice, then this implies
$ (\an) _ k\in {E _ 0}.$
Therefore, by (21),
$(\an) _ k\in {E _ 1}, $
i.e.,
$ (\wx) _ k\in {E _ 1} $
or
$ x\in {X _ 1}. $
Thus
$ X _ 0\subset {X _ 1}. $
Arguing as above, we obtain the converse inclusion, and
$ X _ 0 = X _ 1$
as sets. Since
$X_0$
and
$X_1$
are Banach lattices, then their norms are equivalent.

This completes the proof.
\end{proof}
\vskip 0.5cm

\begin {center}
{\large\bf 3. Final remarks and examples}
\end {center}
\vskip 0.1cm
\begin{rem}
Combining Theorems 1 --- 3 with results obtained in the paper [18], we may
also prove  similar assertions for lacunary trigonometric series. Moreover,
taking into account the main result of the paper [19], we may extend
Theorems 1 --- 3 to Sidon systems of characters of a compact
Abelian group.
\end{rem}

\begin{rem}
In Theorem 1, we cannot replace the space
$ G $
by
$ L _ q $
with some
$ q < \infty. $
Indeed, suppose that the couple
$ (T (l _ 1), T(l _ 2)) $
is a
$ \ab $
-subcouple of the couple
$ (\rv, L _ q),$
i.e.,
$$
\fg \,\asymp {{\ab} (t, Ta; \rv, L _ q)}. $$
Let
$ E = l _ p (2 ^ {-\theta k}),$
where
$ 0 < \theta< 1$
and
$ p = q / {\theta}.$
Applying the
$ \ab $
-method of interpolation
$(\cdot,\cdot)_E^{\ab}$
to the couples
$ (l _ 1, l _ 2) $
and
$ (\rv, L _ q)$
we obtain
$$ || Ta || _ p \,\asymp {\, || a || _ {r, p}} \, = \, \left \{\sum _ {k = 1}^ {\infty} (a _ k ^ *) ^ pk ^ {p/r-1} \right \} ^ {1/p}. $$
Since
$ r = 2 /(2-\theta) < 2$
[4, p.142], then this contradicts with (2).
\end{rem}

\begin{rem}
Clearly, a partial retract of a couple
$ \eu = (Y _ 0, Y _ 1)$
is also a
$ \ab $
-subcouple of
$ \eu.$
The opposite assertion is not true, in general (nevertheless, some
interesting examples of
$ \ab $
-subcouples and partial retracts simultaneously are given in the paper [20]).
Indeed, by Theorem 1, the subcouple
$ (l _ 1, l _ 2) $
is a
$\ab$
-subcouple of the couple
$ (\rv,G). $
Assume that
$ (l _ 1, l _ 2) $
is a partial retract of this couple. Then (see the proof of  Proposition 1)
$ (l _ 1, l _ 2) $
is a partial retract of the couple
$ (\rv, \rv (\log _ 2 ^ {-1/2} (2/t))), $
as well. Therefore, by Lemma 1, from [21] and [4, p.142] it follows that
$$
[ l _ 1, l _ 2] _ {\theta} \, = \, (l _ 1, l _ 2) _ {\theta, \infty} \, = \, l_ {p, \infty}, $$
where
$[ l _ 1, l _ 2] _ {\theta}$
is the space of the complex method of interpolation [4],
$ 0 <\theta < 1,$
and
$ p = \, 2 / (2-\theta). $
On the other hand, it is well-known [4, p.139]
that
$$
[ l _ 1, l _ 2] _ {\theta} \, = \, l _ p \;\;\;\;\mbox{for}\;\;\; p = \, \frac{2}{2-\theta}. $$
This contradiction shows that the couple
$ (l _ 1, l _ 2) $
is not a partial retract of the couple
$ (\rv, G). $
\end{rem}
\vskip 0.15cm

Using Theorem 2, we can find coordinate sequence spaces of
coefficients of Rademacher series from certain r.i.s.'s.

\begin{rem1}
Let $ X $
be the Marcinkiewicz space
$ \xy, $
where
$ \varphi(t) = \, = t\log _ 2\log _ 2 (16/t),$
$0 < t\le 1. $
Show that
$$
\biggl\|\cd\biggr\| _ {\xy} \, \asymp {\, || a || _ {\yz}}, \eqno {(23)} $$
where
$ \yz $
is the space of all sequences
$ a = \gh $
such that
$$ || a || _ {\yz} \, = \, \sup _ {k = 1,2..} \log _ 2 ^ {-1} (2k) \sum _ {i =1} ^ k \, a _ i ^ *\eqno {(24)} $$
is finite. Taking into account Theorem 2, it is sufficient to check that
$$
\ij _ F \, = \,\yz\eqno {(25)} $$
and
$$
{\jk} _ F \, = \,\xy,\eqno {(26)} $$
for some parameter
$F$
of the
$ \ab $
-method of interpolation.
More precisely, we shall prove that (25) and (26) are true for
$ F = l _ {\infty} (u _ k),$
where
$u _ k = 1 / (k + 1)$
($ k\ge 0$)
and
$ u _k = 1 $
($ k < 0).$

By the Holmstedt formula (9),
$$
\kl (2 ^ k) \, \le {\, \sum _ {i = 1} ^ {2 ^ {2k}} a _ i ^ * + 2 ^ k\biggl[\sum _ {i = 2 ^ {2k} + 1} ^ {\infty} (a _ i ^ *) ^ 2\biggr] ^ {1/2}} \, \le{\, B\kl (2 ^ k)}\;\;\;\mbox{for}\;\;\;k=0,1,2,.., \eqno {(27)} $$
where, as before,
$ \kl (t) = \fg.$
Without loss of generality, assume that
$ a _ i = a _ i ^ *. $
If
$ ||a || _ {\yz} = \, R \, <\infty,$
then by (24),
$$
\sum _ {i = 1} ^ {2 ^ {2k}} a _ i ^ * \,\le {\, 2R (k + 1)}. \eqno {(28)} $$
In particular, this implies
$ a _ {2 ^ {2k}} \le {2 ^ {-2k + 1} R (k + 1)}$
for nonnegative integer
$k.$
Using the last inequality, we obtain
\begin{multline*}
\sum _ {i = 2 ^ {2k} + 1} ^ {\infty} a _ i ^ 2 \, = \,\sum _ {j = k} ^{\infty} \sum _ {i = 2 ^ {2j} + 1} ^ {2 ^ {2 (j + 1)}} a _ i ^ 2 \,\le {\,3\sum _ {j = k} ^ {\infty} 2 ^ {2j} a _ {2 ^ {2j}} ^ 2} \, \le \\
\le {12R ^ 2\sum _ {j = k} ^ {\infty} 2 ^ {-2j} (j + 1) ^ 2} \, \le {\, 192R ^2\int _ {k + 1} ^ {\infty} x ^ 22 ^ {-2x} \, dx} \, \le {\, 144R ^ 2 (k + 1) ^2 {2 ^ {-2k}}}.
\end{multline*}
Hence the second term in (27) does not exceed
$ 12R (k + 1). $
Therefore, if
$ E = \, \ij _ F,$
then (28) implies
$$
|| a || _ E \, = \,\sup _ {k = 0,1..} \, \frac {\kl (2 ^ k)} {k + 1} \, \le {\,14 || a || _ {\yz}}. $$

Conversely, if
$ 2 ^ {2j} + 1\le k\le {2 ^ {2 (j + 1)}}$
for some
$j=0,1,2,..,$
then from (27) it follows
$$ \sum _ {i = 1} ^ ka _ i \,\le {\, B\kl (2 ^ {j + 1})} \, \le {\, \sum _ {i =1} ^ {2 ^ {2 (j + 1)}} a _ i} \, \le {\, B || a || _ E (j + 2)} \, \le {\,2B\log _ 2 (2k) || a || _ E.} $$
Therefore
$ || a || _ {\yz} \le {2B || a || _E}$
and (25) is proved.

We pass now to function spaces. At first, we introduce one more
interpolation method which is, actually, a special case of the real
method of interpolation.

For a function
$ \varphi\in {\op}$
and an arbitrary Banach couple
$ (X _ 0, X _ 1)$
define generalized Marcinkiewicz space as follows:
$$ M _ {\varphi} (X _ 0, X _ 1) \, = \, \biggl \{x\in {X _ 0 + X _ 1}: \,\sup _{t> 0} \frac {{\ab} (t, x; X _ 0, X _ 1)} {\varphi (t)}\,<\, \infty\biggr\}. $$
By equation (16), we have
$$
\rv \, = \, M _ {\us} (L _ 1, \rv)\;\;\;\;\mbox{and}\;\;\;\; L _ N \, = \, M _ {\st} (L _ 1, \rv),$$
where these spaces are function spaces on the segment
$[0,1].$
Here
$ \us (t) = \, \min (1, t),$
$\st (t) = \, \min (1, t\log _ 2 ^ {1/2}[\max (2,2/t)]),$
and
$ N (t) = \exp (t ^ 2) -1,$
as before. In addition, using similar notation, it is easy to check that
$$
\dn _ F \, = \, M _ {\rho} (X _ 0, X _ 1),$$
for an arbitrary Banach couple
$ (X _ 0, X _ 1) $
and
$ \rho (t) = \, \log _ 2 (4 + t). $
Hence, by the reiteration theorem for generalized Marcinkiewicz spaces
[11, p.428], we obtain
$$ (L _ {\infty}, L _ N) _ F ^ {\ab} \, = \, M _ {\rho} (M _ {\us} (L _ 1,\rv), M _ {\st} (L _ 1, \rv)) \,=\, M _ {\varphi _ {\rho}} (L _ 1, \rv) \, = \,M (\varphi _ {\rho}),$$
where
$ \varphi _ {\rho} (t) = \, \us (t) \rho (\st (t) /\us(t)). $
A simple calculation gives
$ \varphi _ {\rho} (t) \asymp {\varphi (t)},$
if
$t>0.$
Thus,
$$
( \rv, L _ N) _ F ^ {\ab} \, = \, \xy. $$
It is readily seen that
$ \ef\,=\,{\ab} (t, x; \rv, L _ N), $
for all
$x\in G.$
Therefore, for such
$x$
the norm
$\|x\|_{\xy}$
is equal to the norm
$\|x\|_Y,$
where
$ Y = \jk_ F.$
On the other hand, for
$ x\in {\xy}$
$$ \frac {1} {t\log _ 2 ^ {1/2} (2/t)} \int _ 0 ^ tx ^ * (s) ds \,\le {\, || x|| _ {\xy} \frac {\log _ 2\log _ 2 (16/t)} {\log _ 2 ^ {1/2} (2/t)}} \to 0\;\;\;\mbox{as}\;\;t\to 0+.$$
This implies that
$ \xy\subset G $
[2, p.156].

Thus
$ Y = \xy $
and therefore (26) is proved. Equivalence (23) follows now, as already
stated, from (25) and (26).
\end{rem1}

\begin{rem}
Theorems 2 and 3 strengthen results of the papers [3] and [22],
where similar assertions are obtained for sequence spaces
$ F $
satisfying more restrictive conditions.
For instance, we can readily show that the norm of the dilation operator
$$
\sigma _ na = \, (\underbrace {a _ 1,., a _ 1} _ {n}, \underbrace {a _ 2,., a _ 2} _ {n},..) $$
in the space
$ l _ 1 (\ln) $
(see Example 2) is equal to
$ n $.
Therefore condition (11) from [3] fails for this space
and the theorems obtained in the papers [3] and [22] cannot be applied
to it. Similarly, the Marcinkiewicz space
$ \xy $
from Example 1 does not satisfy conditions of Theorem 8 of [3].
\end{rem}

Using Theorems 2 and 3, we can derive certain interpolation relations.

\begin{rem1}
Let
$\varphi\in {\op}$
and
$1\le p < \infty.$
Recall that the Lorentz space
$ \Lambda _ p(\varphi)$
consists of all measurable functions
$ x = x (s) $
such that
$$
|| x || _ {\varphi, p} \, = \, \left \{\int _ 0 ^ 1 \, (x ^ * (s)) ^ p \, d\varphi(s) \right \} ^ {1/p}\,<\,\infty. $$
In the paper [3], V.A.Rodin and E.M.Semenov proved that
$$
\biggl\|\sum _ {k = 1} ^ {\infty} \, a _ kr _ k \biggr\| _ {\varphi, p} \, \asymp {\, || (a _k) || _ p}, $$
where
$ \varphi (s) = \, \log _ 2 ^ {1-p} (2/s) $
and
$ 1 < p< 2.$
Moreover, the space
$ \Lambda _ p (\varphi) $
is the unique r.i.s. having this property.
Note that
$ l _ p = (l _ 1, l _ 2) _ {\theta, p},$
where
$ \theta = 2 (p-1) /p$
[4, p.142]. Therefore, by Theorem 2,
we obtain
$$
( \rv, G) _ {\theta, p} \, = \, \td,$$
for the same
$ p $
and
$ \theta. $
\end{rem1}
\newpage
\pagestyle{empty}

\end {document}